\documentclass[a4paper, 11pt]{article}

\usepackage{inputenc}
\usepackage{amsmath, amsthm, amssymb}
\usepackage{color}
\usepackage{graphicx}
\usepackage[a4paper]{geometry}

\numberwithin{equation}{section}

\newtheorem{thm}{Theorem}[section]

\theoremstyle{definition}
\newtheorem{expl}[thm]{Example}
\newtheorem{rem}[thm]{Remark}

\DeclareMathOperator{\DD}{\ensuremath\normalfont{D}}

\newcommand{\R}{\mathbb{R}}

\newenvironment{abst}{\begin{minipage}[c]{0.9\textwidth} \footnotesize \textbf{Abstract.}}
{\end{minipage}\\[2ex]}

\newenvironment{key}{\begin{minipage}[c]{0.9\textwidth} \footnotesize \textbf{Keywords.}}
{\end{minipage}\\[2ex]}

\begin{document}
\begin{center}
\Large\bfseries A Parameter Choice Strategy for the Inversion of Multiple Observations \normalsize\mdseries
\\[3ex]
{Christian Gerhards}\footnote{Computational Science Center, University of Vienna, 1090 Vienna, Autria, e-mail: christian.gerhards@univie.ac.at},
{Sergiy Pereverzyev Jr.}\footnote{Institute of Mathematics, University of Innsbruck, 6020 Innsbruck, Austria, e-mail: sergiy.pereverzyev@uibk.ac.at },
{Pavlo Tkachenko}\footnote{ RICAM, 4040 Linz, Austria, e-mail: pavlo.tkachenko@oeaw.ac.at}
\\[3ex]
\end{center}

\begin{abst}

In many geoscientific applications, multiple noisy observations of different origin need to be combined to improve the reconstruction of a common underlying quantity. This naturally leads to multi-parameter models for which adequate strategies are required to choose a set of 'good' parameters. In this study, we present a fairly general method for choosing such a set of parameters, provided that discrete direct, but maybe noisy, measurements of the underlying quantity are included in the observation data, and the inner product of the reconstruction space can be accurately estimated by the inner product of the discretization space. Then the proposed parameter choice method gives an accuracy that only by an absolute constant multiplier differs from the noise level and the accuracy of the best approximant in the reconstruction and in the discretization spaces.
\end{abst}

\begin{key}
Parameter Choice, Multiple Observations, Spherical Approximation
\end{key}

\section{Introduction}

Satellite missions like CHAMP, GRACE, GOCE, or Swarm (e.g., \cite{champ,drinkwater03,friischristensen06,grace}) provide highly accurate data of the Earth's gravity and magnetic field, e.g., by giving information on the first- or second-order radial derivative of the gravitational potential or measurements of the vectorial geomagnetic field, which, once certain iono- and magnetospheric contributions have been filtered out, can be expressed as the gradient of a harmonic potential. Drawing conclusions from such satellite measurements on the gravitational potential or the magnetic field at or near the Earth's surface is a classical exponentially ill-posed problem (see, e.g., \cite{freeden99,lu15,pereverzyev99}). Measurements at or near the Earth's surface (which we simply denote as ground measurements), on the other hand, do not suffer from this ill-posedness but are typically only available in restricted regions (e.g., aeromagnetic surveying \cite{pikington07}). Combining both data sets becomes necessary when aiming at local high resolution models that also take global trends into account. This is a classical setting for multiparameter modeling (e.g., \cite{brezinski03,koch02,lu11,lu15})  that involves the regularization of an ill-posed inverse problem (downward continuation of satellite data) and the weighting of the satellite data against the ground data. An exemplary situation that we also use for later numerical illustrations is the following: We have measurements $f_1$ of a harmonic potential $u$ on a spherical satellite orbit $\Omega_R=\{x\in\mathbb{R}^3:|x|=R\}$ and measurements $f_2$ of $u$ in a subregion $\Gamma_r\subset\Omega_r$ of the spherical Earth's surface $\Omega_r$, $r<R$, i.e.,
\begin{align}
 \Delta u&=0, \textnormal{ in }\Omega_r^{ext},\label{eqn:mag1}
 \\u&=f_2, \textnormal{ on }\Omega_R,
 \\u&=f_1,\textnormal{ on }\Gamma_r,\label{eqn:mag4}
\end{align}
where $\Omega_r^{ext}=\{x\in\mathbb{R}^3:|x|>r\}$. The problem of approximating $u$ in $\Gamma_r$ is clearly overdetermined and spherical splines (e.g., \cite{freeden81,shure82}) or other localized basis functions (e.g., \cite{haines85,simons06,thebault06}) could be used to approximate $u$ in $\Gamma_r$ from knowledge of $f_1$ only (generally, we denote the restriction of $u$ to $\Gamma_r$ by $u^\dag$). However, such methods are not always well-suited to capture global trends of $u$ and they do not address situations where the noise level of $f_2$ might be smaller than that of $f_1$. Therefore, it is advisable to incorporate satellite data $f_2$ as well. Eventually, based on different parameter settings or approximation methods, we assume to have a set of candidates $\{u_k\}_{k=1,2,\ldots,N}$ available for the approximation of $u$ in $\Gamma_r$. 

In this paper, we aim at introducing a method that predicts a 'good' candidate $u_{k^*}$ among the available $\{u_k\}_{k=1,2,\ldots,N}$ without requiring knowledge of the method by which each $u_k$ has been obtained or which sort of noise is contained in the data. It is also not necessary to know the underlying models or the type of data that has lead to the construction of $u_k$. Apart from $\{u_k\}_{k=1,2,\ldots,N}$, all that is required is a reference measurement $f$ (in the example \eqref{eqn:mag1}--\eqref{eqn:mag4}, this would be $f_1$) of $u^\dag$ against which to compare the candidates $u_k$. In this sense, we are not dealing with a parameter choice strategy for an ill-posed problem (although the underlying models that determine $u$ may be ill-posed) but rather with a general method of choosing a 'good' approximant of $u^\dag$ among a set of available candidates
(an extensive comparison of parameter choice methods for ill-posed problems can be found, e.g., in \cite{bauer15,bauer11}).
Opposed to aggregation methods (see, e.g., \cite{chen15,merrick13}), where approximations from different data settings are superposed to obtain a final approximation, we assume in our method that this superposition has already taken place in one way or another during the construction of each $u_k$. An important constraint for our method, in order to obtain a suitable error estimate, is that the discrete reference measurements $f$ of $u^\dag$ need to be given in points that allow the definition of an inner product in the discretization space which coincides with the $L^2$-inner product in a desired finite-dimensional function space (e.g., the spherical harmonics of degree smaller than some $L$). The numerical tests, however, show that our method also supplies good results if this condition is slightly violated.

The structure of the paper is as follows: In Section \ref{sec:theory}, we introduce and investigate the parameter choice strategy mentioned above in more detail and put it into a mathematically rigorous context. In Section \ref{sec:numerics}, we illustrate its performance for the problem \eqref{eqn:mag1}--\eqref{eqn:mag4}. The approximations for this problem are obtained by a method described in \cite{gerhards14b}. Latter is also briefly recapitulated in Section \ref{sec:numerics}.

\section{The Parameter Choice Strategy}\label{sec:theory}

Throughout this paper, we assume the following conditions to be satisfied:
\begin{itemize}
\item[(a)] Let $\Gamma_r\subset\Omega_r$ be a subdomain of the sphere $\Omega_r$, where discrete direct measurements of the underlying quantity $u$ are available. We assume to have $M$ measurement values
and a corresponding discretization operator $\DD:L^2(\Gamma_r)\rightarrow \mathbb{R}^M$, that maps a function $u^\dag\in L^2(\Gamma_r)$ to the correponding measurements ${\DD}u^\dag\in\mathbb{R}^M$. Furthermore, let $\mathbb{R}^M$ be equipped with some inner product $\left\langle \cdot, \cdot\right\rangle_{\mathbb{R}^M}$ and the corresponding norm $\left\|\cdot\right\|_{\mathbb{R}^M}$.

\item[(b)] The measurements of $u^\dag$ may be blurred by additive noise $\xi=(\xi_1,\ldots,\xi_M)\in\R^M$ and we assume, without loss of generality, that there is $u_\xi^\dag =f \in L^2(\Gamma_r)$ such that
\begin{align*}
{\DD}u_\xi^\dag={\DD}u^\dag+\xi,\quad \left\|{\DD}u^\dag-{\DD}u_\xi^\dag\right\|_{\R^M}\leq \varepsilon
\end{align*}
for some $\varepsilon>0$.

\item[(c)] We assume that from somewhere, a set $\left\{u_k\right\}_{k=1,2,\ldots,N}$ of approximations of $u^\dag$ on $\Gamma_r$ is available sand that all these approximations belong to some finite dimensional linear subspace $V\subset L^2(\Gamma_r)$.

\item[(d)] Finally, we assume that the discretization space $\R^M$ is related to the reconstruction space $V$ through the discretization operator $\DD$ such that
\begin{align}\label{cond_3}
\left\langle g, \bar{g}\right\rangle_{L^2(\Gamma_r)}=\left\langle {\DD}g, {\DD}\bar{g}\right\rangle_{\R^M},\quad\textnormal{for all }g, \bar{g} \in V.
\end{align}
\end{itemize}

\begin{expl}\label{example}
Let $V=V_L$ be the space of spherical polynomials of the degree $L$. Under rather general assumptions on $\Gamma_r$ one can find a system of knots $\left\{x_i^M\right\}_{i=1,\ldots,M}$ and positive weights $\left\{w_i^M\right\}_{i=1,\ldots,M}$ such that
\begin{align*}
\int_{\Gamma_r}g(x)d\Gamma_r(x)=\sum_{i=1}^M w_i^Mg(x_i^M),\quad \textnormal{for all }g\in V_{2L}.
\end{align*}
Consider a discretization operator
\begin{align*}
{\DD}g=(g(x_1^M), g(x_2^M),\ldots, g(x_M^M))\in\R^M
\end{align*}
and the inner product
\begin{align*}
\left\langle y, \bar{y}\right\rangle_{\R^M}:=\sum_{i=1}^M w_i^M y_i\bar{y}_i.
\end{align*}
It is clear that for the just introduced reconstruction space $V=V_L$, discretization space $\R^M$, and discretization operator $\DD$ the condition \eqref{cond_3} is satisfied. It is also clear that the measurements described by the operator $\DD$ are just pointwise evaluations at the knots $\left\{x_i^M\right\}_{i=1,\ldots,M}$, and the noise level of these measurements is controlled by the quantity
\begin{align*}
\left\|{\DD}u^\dag-{\DD}u_\xi^\dag\right\|_{\R^M}^2=\sum_{i=1}^M w_i^M (u^\dag(x_i^M)-u_\xi^\dag(x_i^M))^2\leq \varepsilon^2.
\end{align*}
\end{expl}

Now we are ready to describe our choice of a good approximation from a given family $\left\{u_k\right\}_{k=1,\ldots,N}$. Let $u_{k_{opt}}\in\left\{u_k\right\}_{k=1,\ldots,N}$ be such that
\begin{align*}
\left\|u^\dag-u_{k_{opt}}\right\|_{L^2(\Gamma_r)}=\min_{k=1,2,\ldots,N}\left\|u^\dag-u_k\right\|_{L^2(\Gamma_r)}.
\end{align*}
Of course, $u_{k_{opt}}$ cannot be found without knowledge of $u^\dag$. Therefore, in practice one cannot find the parameters corresponding to  $u_{k_{opt}}$.

We motivate our procedure with the observation that for any $u_k$, $k=1,2,\ldots,N$, it holds
\begin{align*}
	\left\|u_k-u_{k_{opt}}\right\|_{L^2(\Gamma_r)}&=\sup_{a\in L^2(\Gamma_r),\|a\|_{L^2(\Gamma_r)}=1}\left\langle u_k-u_{k_{opt}},a\right\rangle_{L^2(\Gamma_r)}
	\\&=
	\max_{a\in A_N}\left\langle u_k-u_{k_{opt}},a\right\rangle_{L^2(\Gamma_r)},
\end{align*} 
where the finite set $A_N$ is defined as follows
\begin{align*}
A_N=\left\{a=a_{k,l}=\frac{u_k-u_l}{\left\|u_k-u_l\right\|_{L^2(\Gamma_r)}},\ k, l=1,2,\ldots,N\right\}\subset V.
\end{align*}
Then for any $k=1,2,\ldots,N$ and $a\in A_N$ the quantity 
\begin{align*}
\left\langle u_k-u_{k_{opt}},a\right\rangle_{L^2(\Gamma_r)}=  \left\langle u_k,a\right\rangle_{L^2(\Gamma_r)}-\left\langle u_{k_{opt}},a\right\rangle_{L^2(\Gamma_r)}
\end{align*}
has only a part $\left\langle u_{k_{opt}},a\right\rangle_{L^2(\Gamma_r)}$ that cannot be computed directly, because $u_{k_{opt}}$ is unknown. On the other hand, this part can be approximated with the use of the available observations as follows
\begin{align*}
\left\langle u_{k_{opt}},a\right\rangle_{L^2(\Gamma_r)}=\left\langle {\DD}u_{k_{opt}},{\DD}a\right\rangle_{\R^M}\approx \left\langle {\DD}u_\xi^\dag,{\DD}a\right\rangle_{\R^M}.
\end{align*}
Therefore, the values
\begin{align*}
h_k(a)=\left\langle u_k,a\right\rangle_{L^2(\Gamma_r)}-\left\langle {\DD}u_\xi^\dag,{\DD}a\right\rangle_{\R^M}
\end{align*}
and
\begin{align*}
H_k=\max_{a\in A_N}\left|h_k(a)\right|
\end{align*}
can be seen as surrogates for the values of $\left\langle u_k-u_{k_{opt}},a\right\rangle_{L^2(\Gamma_r)}$ and $ \left\|u_k-u_{k_{opt}}\right\|_{L^2(\Gamma_r)}$ respectively.

In the view of this it is natural to expect that the approximation $u_{k_*}\in \left\{u_k\right\}_{k=1,\ldots,N}$ defined by
\begin{align}\label{oracle_index}
k_*: \ H_{k_*}=\min\left\{H_k,\ k=1,2,\ldots,N\right\}
\end{align}
is close to $u_{k_{opt}}$. Indeed,
\begin{align} \label{bound}
\left\|u_{k_*}-u_{k_{opt}}\right\|_{L^2(\Gamma_r)}&=\left\langle u_{k_*}-u_{k_{opt}},a_{k_*,k_{opt}}\right\rangle_{L^2(\Gamma_r)} \\ \nonumber
&=\left(\left\langle u_{k_*},a_{k_*,k_{opt}}\right\rangle_{L^2(\Gamma_r)}-\left\langle {\DD}u_\xi^\dag,a_{k_*,k_{opt}}\right\rangle_{\R^M}\right) \\ \nonumber
&-\left(\left\langle u_{k_{opt}},a_{k_*,k_{opt}}\right\rangle_{L^2(\Gamma_r)}-\left\langle {\DD}u_\xi^\dag,a_{k_*,k_{opt}}\right\rangle_{\R^M}\right) \\ \nonumber
&=h_{k_*}(a_{k_*,k_{opt}})-h_{k_{opt}}(a_{k_*,k_{opt}}) \\ \nonumber
&\leq H_{k_*}+H_{k_{opt}}\leq 2H_{k_{opt}}.
\end{align}
Furthermore,
\begin{align*}
H_{k_{opt}}&=\max_{a\in A_N} \left|\left\langle u_{k_{opt}},a\right\rangle_{L^2(\Gamma_r)}-\left\langle {\DD}u_\xi^\dag,{\DD}a\right\rangle_{\R^M}\right| \\
&=\max_{a\in A_N} \left|\left\langle {\DD}u_{k_{opt}}-{\DD}u^\dag,{\DD}a\right\rangle_{\R^M}+\left\langle {\DD}u^\dag-{\DD}u_\xi^\dag,{\DD}a\right\rangle_{\R^M}\right| \\
&\leq \varepsilon+\left\|{\DD}u^\dag-{\DD}u_{k_{opt}}\right\|_{\R^M}.
\end{align*}

Then the previous analysis together with the triangle inequality gives us the following statement.

\begin{thm}\label{thm}
Let us assume that conditions (a)--(d) hold true, i.e., we are given a finite family of approximations $\left\{u_k\right\}_{k=1,\ldots,N}$ from a finite dimensional reconstruction space $V\subset L^2(\Gamma_r)$. Moreover, noisy direct discrete measurements ${\DD}u_\xi^\dag\in \R^M$ of the approximated quantity $u^\dag$ are available, and the reconstruction space $V$ is related to the discretization space $\R^M$ such that \eqref{cond_3} is satisfied. Then for $k_*$ chosen according to (\ref{oracle_index}) we have
\begin{align*}
\left\|u^\dag-u_{k_*}\right\|_{L^2(\Gamma_r)}\leq \left\|u^\dag-u_{k_{opt}}\right\|_{L^2(\Gamma_r)}+2\left\|{\DD}u^\dag-{\DD}u_{k_{opt}}\right\|_{\R^M}+2\varepsilon.
\end{align*}
\end{thm}

\begin{rem}\label{rem}
Note that in the context of Example \ref{example} we can give also another bound for $\left\|u^\dag-u_{k_*}\right\|_{L^2(\Gamma_r)}$. Let $u^L_{best}\in V_L$ be the spherical polynomial of the best $C(\Gamma_r)$-approximation, i.e.,
\begin{align*}
\left\|u^\dag-u^L_{best}\right\|_{C(\Gamma_r)}=\min_{v\in V_L}\left\|u^\dag-v\right\|_{C(\Gamma_r)}.
\end{align*}
Then
\begin{align*}
h_{k_{opt}}(a)&=\left\langle u^\dag, a\right\rangle_{L^2(\Gamma_r)}+\left\langle u_{k_{opt}}-u^\dag,a\right\rangle_{L^2(\Gamma_r)}-\left\langle {\DD}u^\dag, {\DD}a\right\rangle_{\R^M} \\
&+\left\langle {\DD}u^\dag-{\DD}u_\xi^\dag, {\DD}a\right\rangle_{\R^M}
\\&=\left\langle u^\dag-u^L_{best}, a\right\rangle_{L^2(\Gamma_r)}-\left\langle {\DD}u^\dag-{\DD}u_{best}^L, {\DD}a\right\rangle_{\R^M} \\
&+\left\langle {\DD}u^\dag-{\DD}u_\xi^\dag, {\DD}a\right\rangle_{\R^M}+\left\langle u_{k_{opt}}-u^\dag,a\right\rangle_{L^2(\Gamma_r)} \\
&\leq\varepsilon+\left\|u^\dag-u_{k_{opt}}\right\|_{L^2(\Gamma_r)}+c_{M,N}\left\|u^\dag-u^L_{best}\right\|_{C(\Gamma_r)},
\end{align*}
where
\begin{align*}
c_{M,N}=\sqrt{\textnormal{vol}(\Gamma_r)}+\sum_{i=1}^M w_i^M\left|a(x_i^M)\right|.
\end{align*}
Furthermore, we observe that
\begin{align*}
\sum_{i=1}^M w_i^M\left|a(x_i^M)\right|\leq \left(\sum_{i=1}^M w_i^Ma^2(x_i^M)\right)^{1/2}\left(\sum_{i=1}^M w_i^M\right)^{1/2}=\sqrt{\textnormal{vol}(\Gamma_r)}.
\end{align*}
Thus,
\begin{align*}
H_{k_{opt}}\leq \varepsilon+\left\|u^\dag-u_{k_{opt}}\right\|_{L^2(\Gamma_r)}+2\sqrt{\textnormal{vol}(\Gamma_r)}\left\|u^\dag-u^L_{best}\right\|_{C(\Gamma_r)},
\end{align*}
and from (\ref{bound}) we get the following alternative bound:
\begin{align*}
\left\|u^\dag-u_{k_*}\right\|\leq 3\left\|u^\dag-u_{k_{opt}}\right\|_{L^2(\Gamma_r)}+4\sqrt{\textnormal{vol}(\Gamma_r)}\left\|u^\dag-u^L_{best}\right\|_{C(\Gamma_r)}+2\varepsilon.
\end{align*}

\end{rem}

\begin{rem}
For the estimates in Theorem \ref{thm} and Remark \ref{rem}, the assumption \eqref{cond_3} is crucial. Incorporating the worst-case error for non-exact quadrature rules in the reconstruction space $V$, estimates similar to those in Theorem \ref{thm} and Remark \ref{rem} could be derived even if \eqref{cond_3} is violated (compare, e.g., \cite{hesse10} for an overview on spherical quadrature rules). However, such estimates would show a stronger and undesirable dependence on $u^\dag$ and $u_k$ and are therefore omitted. In the numerical examples in the next section, we show that a 'slight' violation of the condition \ref{cond_3} can still yield good results.
\end{rem}

\section{Numerical Illustrations}\label{sec:numerics}
In this section, we illustrate the numerical performance of the previously described parameter choice method on \eqref{eqn:mag1}--\eqref{eqn:mag4}, i.e., we assume $u$ to satisfy
\begin{align*}
\Delta u&=0, \textnormal{ in }\Omega_r^{ext},
 \\u&=f_2, \textnormal{ on }\Omega_R,
 \\u&=f_1,\textnormal{ on }\Gamma_r.
\end{align*}
A set of approximations $u_k$, $k=1,\ldots,N$, of $u^\dag$ on a spherical cap $\Gamma_r=\Gamma_r^\rho=\{x\in\Omega_r:1-\frac{x}{|x|}\cdot(0,0,1)^T<\rho\}$ of radius $\rho\in(0,2)$ around the North Pole $(0,0,r)^T$ can be obtained by
\begin{align}\label{eqn:uk}
 u_k(x)=\int_{\Omega_R}\Phi_k(x,y)f_2(y)d\Omega_R(y)+\int_{\Gamma_r}\tilde{\Psi}_k(x,y)f_1(y)d\Gamma_r(y),\quad x\in\Gamma_r,
\end{align}
where the kernels $\Phi_k$, $\Psi_k$ are given by
\begin{align*}
\Phi_k(x,y)&=\sum_{n=0}^{N_k}\sum_{j=1}^{2n+1}\Phi_k^\wedge(n)\frac{1}{r}Y_{n,j}\left(\frac{x}{|x|}\right)
																															\frac{1}{R}Y_{n,j}\left(\frac{y}{|y|}\right),
\\
\tilde{\Psi}_k(x,y)&=\sum_{n=0}^{M_k}\sum_{j=1}^{2n+1}\tilde{\Psi}_k^\wedge(n)
																							\frac{1}{r}Y_{n,j}\left(\frac{x}{|x|}\right)
																							\frac{1}{r}Y_{n,j}\left(\frac{y}{|y|}\right),
\end{align*}
with coefficients $\tilde{\Psi}_k^\wedge(n)$ of the form $\tilde{\Psi}_k^\wedge(n)=\tilde{\Phi}_k^\wedge(n)-\Phi_k^\wedge(n)\left(\frac{r}{R}\right)^n$ and $N_k\leq M_k$.
By $\{Y_{n,j}\}_{n=0,1,\ldots; j=0,1,\ldots,2n+1}$ we mean a set of orthonormal spherical harmonics of degree $n$ and order $j$. The coefficients ${\Phi}_k^\wedge(n)$, $\tilde{\Phi}_k^\wedge(n)$ are chosen by minimizing the following functional:
\begin{equation*}
{\begin{aligned}\label{eqn:mineq2}
\mathcal{F}(\Phi_k,\tilde{\Psi}_k)=&\tilde{\alpha}_{k}\sum_{n=0}^{M_k}\big|1-\tilde{\Phi}_{k}^\wedge(n)\big|^2+\alpha_{k}\sum_{n=0}^{N_k}\left|1-{\Phi}_{k}^\wedge(n)\left(\frac{r}{R}\right)^n\right|^2
\\&+\beta_k\sum_{n=0}^{N_k}\big|\Phi_k^\wedge(n)\big|^2+\big\|\tilde{\Psi}_k\big\|^2_{L^2(\Omega_r\setminus\Gamma_r)}.
\end{aligned}}
\end{equation*}
The first two terms of the functional $\mathcal{F}$ measure the approximation property of the kernels $\Phi_k$, $\tilde{\Psi}_k$ (i.e., they measure how close they are to the Dirichlet kernel). The third term penalizes the error amplification due to the downward continuation of the satellite data on $\Omega_R$ while the fourth term penalizes the localization of $\tilde{\Psi}_k$ outside the region $\Gamma_r$ where ground data is available. The parameters $\alpha_k$, $\tilde{\alpha}_k$, $\beta_k$ weigh these quantities against each other. For more details, on this approach of approximating $u$ on $\Gamma_r$, the reader is referred to \cite{gerhards14b}. Essentially, we are in the setting of Example~\ref{example}
where the reconstruction space $V=V_{M}$ is the space of all spherical polynomials up to degree
$M=\max\left\{ M_k,\; k=1,\ldots,N \right\}$.

\begin{figure}
 \scalebox{0.43}{\includegraphics*{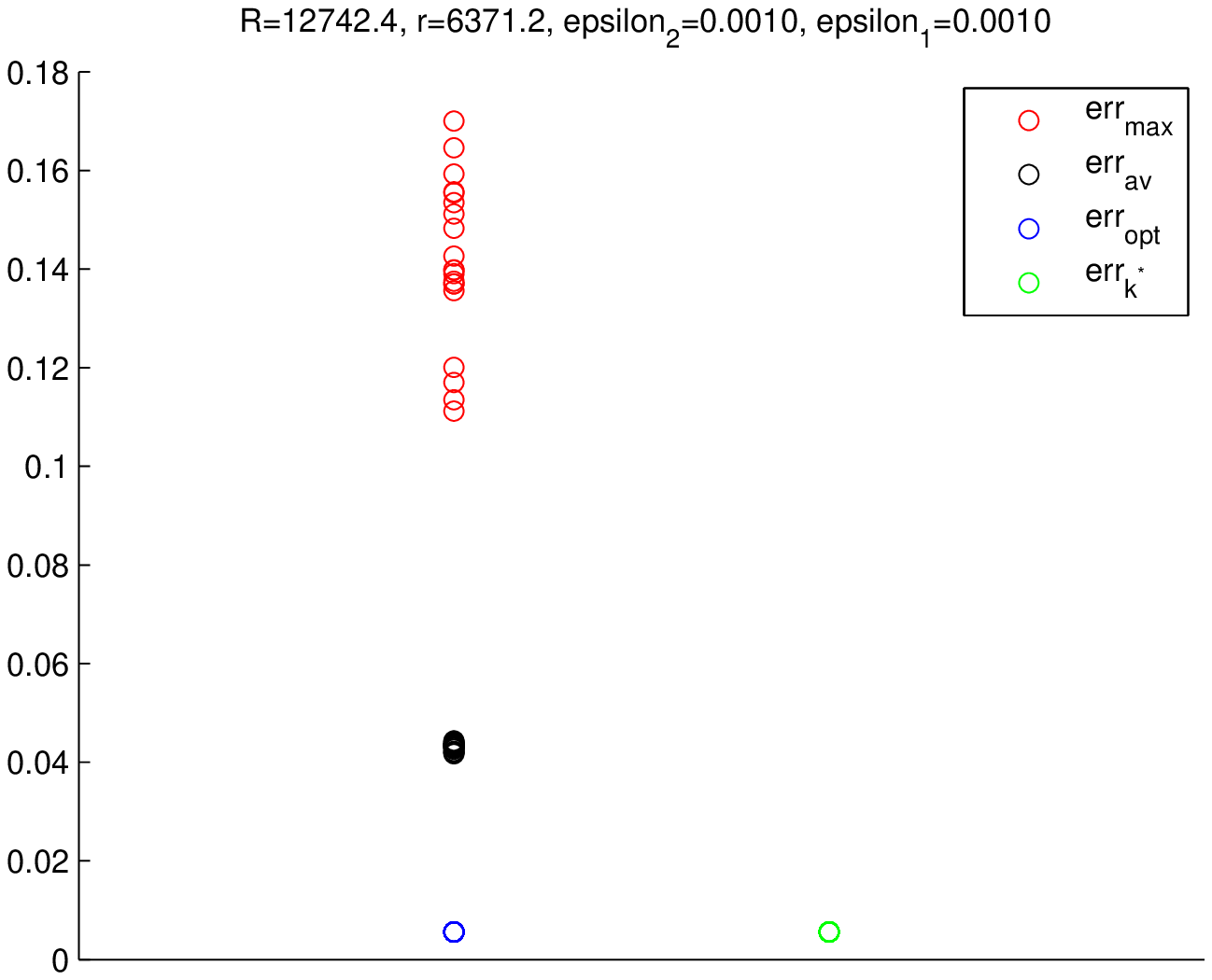}} \scalebox{0.43}{\includegraphics*{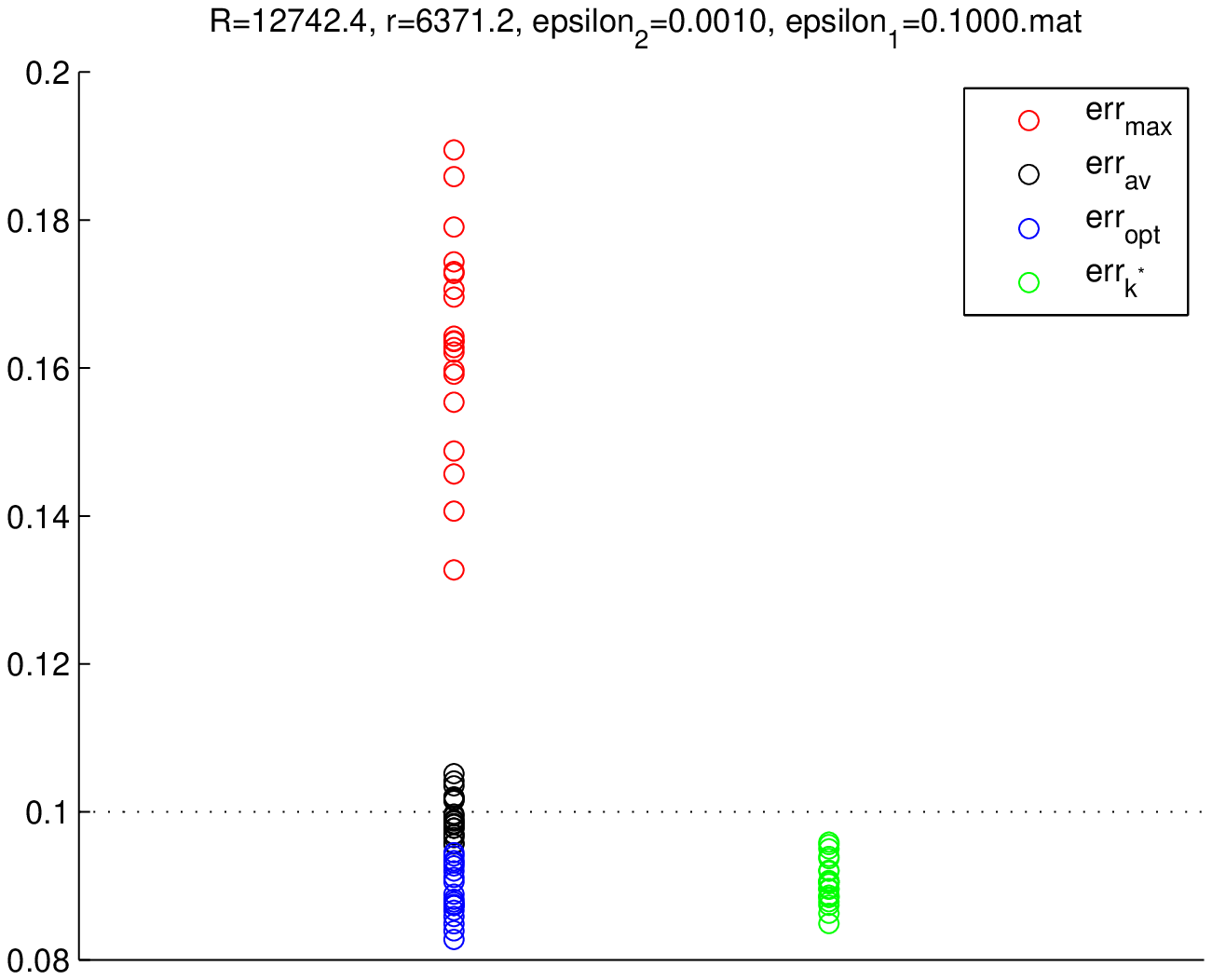}} 
\caption{Relative Errors for the low spherical harmonic degree tests (Situation (a)(1)) for $\varepsilon_1=\varepsilon_2=0.001$ (left) and $\varepsilon_1=0.1$, $\varepsilon_2=0.001$ (right; the dotted black line marks the noise level $\varepsilon_1=0.1$).}\label{fig:test1}
 \scalebox{0.43}{\includegraphics*{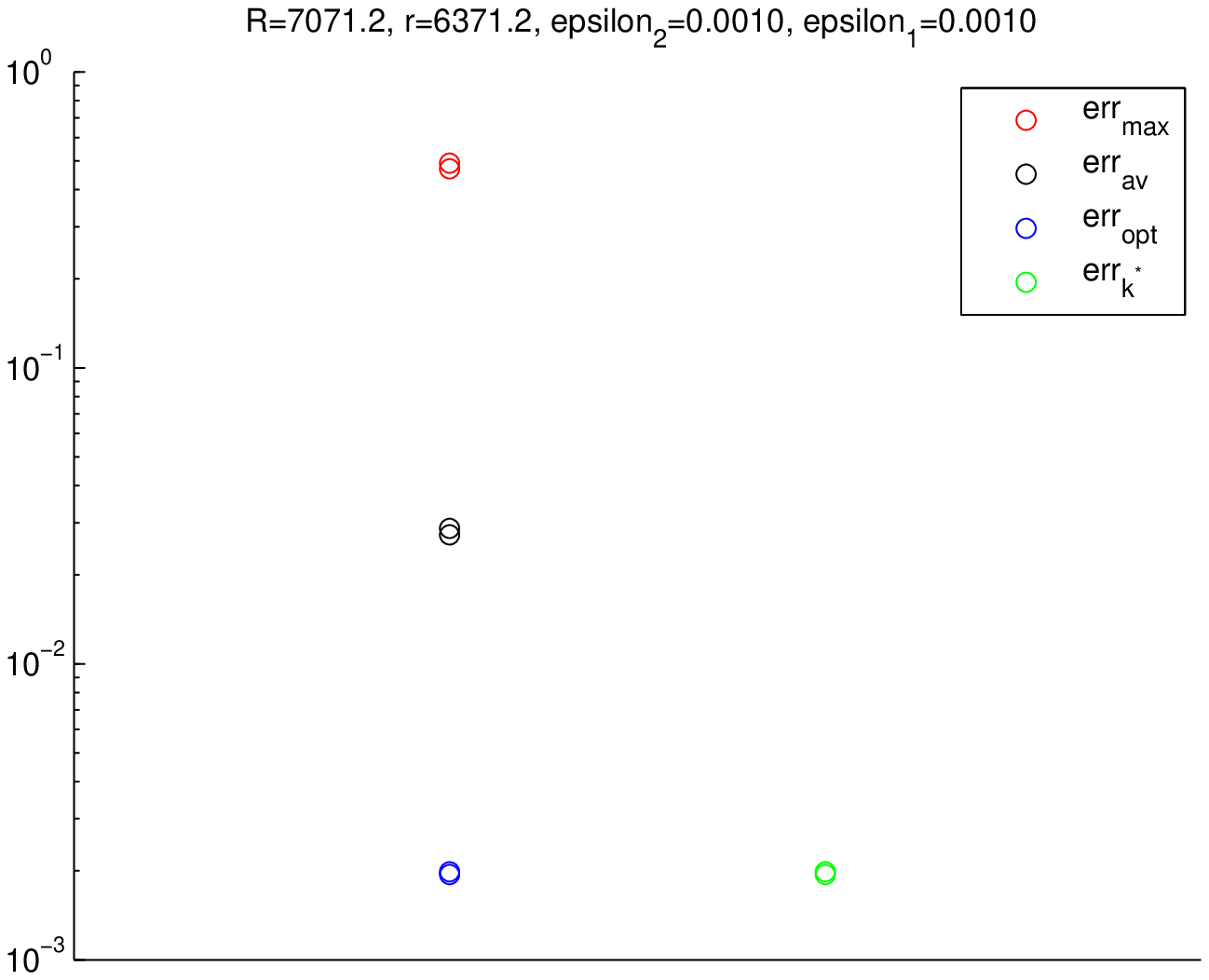}} \scalebox{0.43}{\includegraphics*{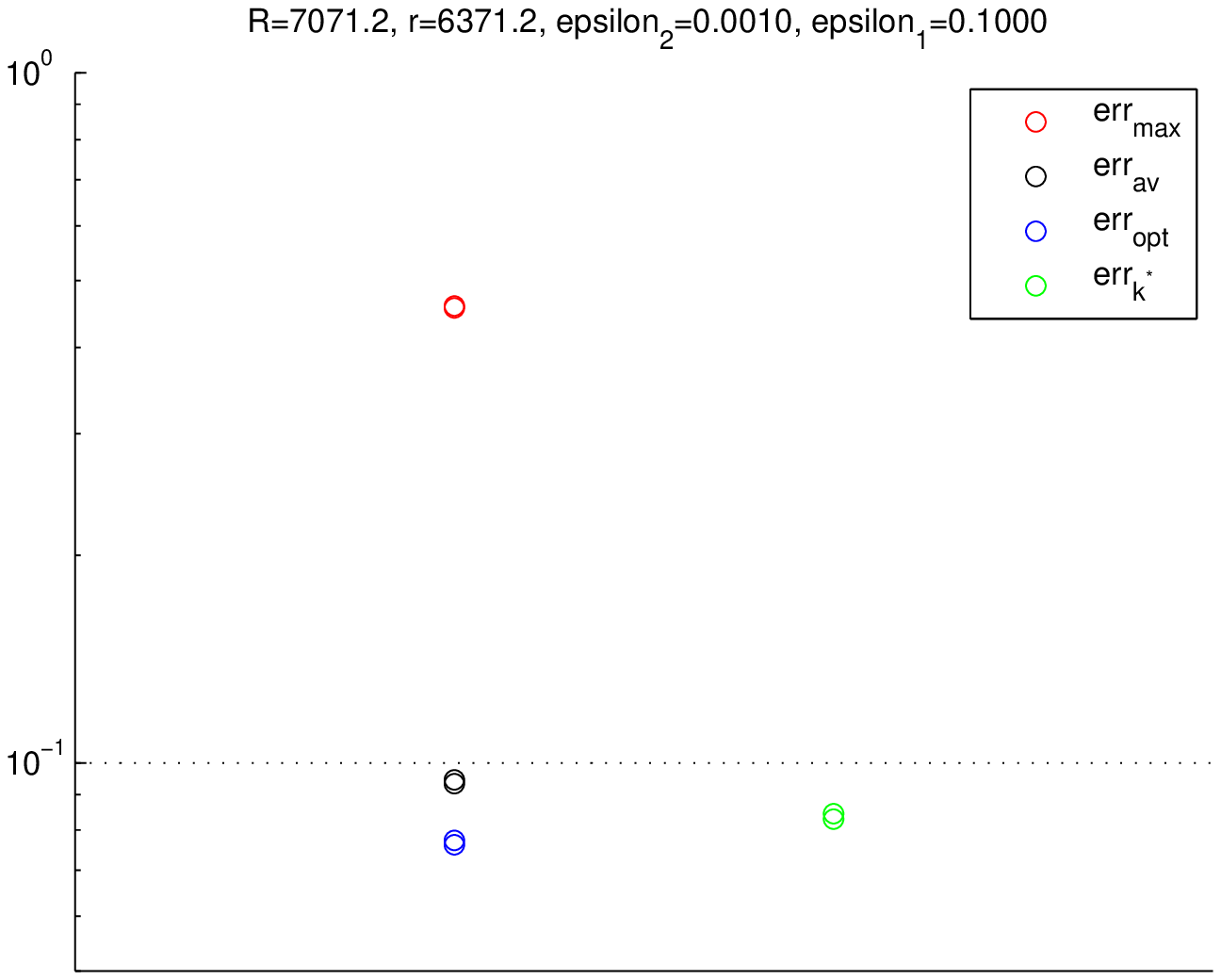}} 
 \caption{Relative Errors for the high spherical harmonic degree tests (Situation (a)(2)) for $\varepsilon_1=\varepsilon_2=0.001$ (left) and $\varepsilon_1=0.1$, $\varepsilon_2=0.001$ (right; the dotted black line marks the noise level $\varepsilon_1=0.1$).}\label{fig:test2}
\end{figure}

The procedure for our numerical tests is as follows:
\begin{itemize}
 \item[(a)] From the EGM2008 gravity potential model (cf. \cite{egm2008}\footnote{data accessed via http://earth-info.nga.mil/GandG/wgs84/gravitymod/egm2008/egm08\_wgs84.html}), we generated two sets of reference potentials $u$: 
 \begin{enumerate}
  \item[(1)] one up to spherical harmonic degree $n=30$ (in order to allow many test runs in a short time) on a sphere $\Omega_R$, $R=12,371$km, and on a spherical cap $\Gamma_r=\Gamma_r^\rho$, $r=6,371$km, with $\rho=1$ (corresponding to a spherical radius of approximately $10,000$km at the Earth's surface),
  \item[(2)] another one up to spherical harmonic degree $n=130$ (in order to have a more realistic scenario) on a sphere $\Omega_R$, $R=7,071$km, and on a spherical cap $\Gamma_r=\Gamma_r^\rho$, $r=6,371$km, with $\rho=0.3$ (corresponding to a spherical radius of approximately $5,000$km at the Earth's surface).
  \end{enumerate}
 \item[(b)] For both cases of part (a), we generate corresponding noisy measurements $f_1$, $f_2$, where the noise levels $\varepsilon_1={\|f_1-u\|_{L^2(\Gamma_r)}}/{\|u\|_{L^2(\Gamma_r)}}$ of the ground data and $\varepsilon_2={\|f_2-u\|_{L^2(\Omega_R)}}/{\|u\|_{L^2(\Omega_R)}}$ of the satellite data are varied among $0.001$, $0.1$. The data on $\Omega_R$ are in both cases computed on an equiangular grid according to \cite{driscoll94} while the data on the spherical cap $\Gamma_r$ are computed on a Gauss-Legendre grid according to \cite{hesse12} in order to guarantee polynomially exact quadrature rules up to spherical polynomial degree $M_k+n$, where $n=30$ in case (1) and $n=130$ in case (2), which yields condition \eqref{cond_3}.
 \item[(c)] For the different input data from part (b), we compute approximations $u_k$, $k=1,\ldots,N$, of $u^\dag$ on $\Gamma_r$ via the expression \eqref{eqn:uk}. The index $k$ of $u_k$ indicates different choices of the parameters $\alpha_k$, $\tilde{\alpha}_k$, $\beta_k$ in the functional $\mathcal{F}$ from \eqref{eqn:mineq2}. $\alpha_k$, $\tilde{\alpha}_k$ are varied in the interval $[10^{1},10^8]$ and $\beta_k$ is varied in the interval $[10^{-2},10^{3}]$. The truncation degrees of the series expansions of $\Phi_k$, $\tilde{\Psi}_k$ are fixed to $N_k=M_k=80$ in case (1) while $N_k=M_k=150$ in case (2). This way, we obtain $N=100$ different approximations $u_k$ for each of the two cases.
 \item[(d)] Among the approximations $u_k$, we use the procedure from Section \ref{sec:theory} to choose a 'good' approximation $u_{k^*}$. Afterwards, we compare the relative approximation errors err$_{k^*}={\|u_{k^*}-u\|_{L^2(\Gamma_r)}}/{\|u\|_{L^2(\Gamma_r)}}$ of the parameter choice with the relative errors err$_{\textnormal{opt}}=\min_{k=1,\ldots,N}{\|u_{k}-u\|_{L^2(\Gamma_r)}}/{\|u\|_{L^2(\Gamma_r)}}$ of the actually best $u_{k_{opt}}$. 
\end{itemize}
The results of the tests are shown in Figures \ref{fig:test1} and \ref{fig:test2}. Each figure shows the relative errors err$_{k^*}$ and err$_{\textnormal{opt}}$ for every test run. Additionally, we plotted the maximum errors err$_{\textnormal{max}}=\max_{k=1,\ldots,N}{\|u_{k}-u\|_{L^2(\Gamma_r)}}/{\|u\|_{L^2(\Gamma_r)}}$ and the average errors err$_{\textnormal{av}}=\frac{1}{N}\sum_{k=1,\ldots,N}{\|u_{k}-u\|_{L^2(\Gamma_r)}}/{\|u\|_{L^2(\Gamma_r)}}$ in order to illustrate the performance. It can be seen that the algorithm works particularly well for the setting $\varepsilon_1=\varepsilon_2$ and that the oracle error err$_{k^*}$ is nearly identical with the minimum error err$_{\textnormal{opt}}$. The situation is different when $\varepsilon_1\gg \varepsilon_2$. The minimum error err$_{\textnormal{opt}}$ is smaller than the noise level $\varepsilon_1$. Thus, since our parameter choice strategy is based on comparing $u_k$ to $f_1$, we cannot expect that err$_{k^*}$ is as good as err$_{\textnormal{opt}}$. Yet, astonishingly enough, it seems that err$_{k^*}$ is still slightly smaller than $\varepsilon_1$ for our test setting.

In addition, we repeated the tests above with a reduced accuracy of the quadrature rule in order to illustrate the consequences if condition \eqref{cond_3} is not satisfied. More precisely, we did test runs for the following setting:
\begin{itemize}
 \item[(a')] We generated two sets of reference potentials $u$: 
 \begin{enumerate}
  \item[(1)] one up to spherical harmonic degree $n=30$ on a sphere $\Omega_R$, $R=12,371$km, and on a spherical cap $\Gamma_r=\Gamma_r^\rho$, $r=6,371$km, with $\rho=1$ (opposed to the previous tests, the potential is not based on the EGM2008 model but the Foruier coefficients are chosen randomly),
  \item[(2)] another one up to spherical harmonic degree $n=130$ on a sphere $\Omega_R$, $R=7,071$km, and on a spherical cap $\Gamma_r=\Gamma_r^\rho$, $r=6,371$km, with $\rho=0.3$ (here, the potential is again based on the EGM2008 model).
  \end{enumerate}
 \item[(b')] For both cases of part (a'), we generate corresponding noisy measurements $f_1$, $f_2$ with noise levels $\varepsilon_1=\varepsilon_2=0.001$. Again, the data on $\Omega_R$ are in both cases computed on an equiangular grid according to \cite{driscoll94} while the data on the spherical cap $\Gamma_r$ are computed on a Gauss-Legendre grid according to \cite{hesse12}. For case (1), we chose grids of two different sizes: one such that the polynomial exactness of the quadrature rule is of degree 100 and one such that polynomial exactness is of degree 90 (remember that polynomial exactness up to degree $M_k+n=80+30=110$ is required in order to satisfy condition \eqref{cond_3}). For case (2), we chose the size of the grids such that the polynomial exactness of the quadrature rule is of degree 130 and of degree 80, respectively (remember that polynomial exactness up to degree $M_k+n=150+130=280$ is required in order to satisfy condition \eqref{cond_3}).
 \item[(c')] For the different input data from part (b), we compute approximations $u_k$, $k=1,\ldots,N$, of $u^\dag$ on $\Gamma_r$ via the expression \eqref{eqn:uk}. The parameters $\alpha_k$, $\tilde{\alpha}_k$ are again varied in the interval $[10^{1},10^8]$ and $\beta_k$ is varied in the interval $[10^{-2},10^{3}]$. The truncation degrees of the series expansions of $\Phi_k$, $\tilde{\Psi}_k$ are fixed to $N_k=M_k=80$ in case (1) while $N_k=M_k=150$ in case (2).  
\end{itemize}
The results are shown in Figures \ref{fig:test3} and \ref{fig:test4}. In the right plot of Figure \ref{fig:test3} it becomes clear that a too large deviation of the required polynomial exactness can severely influence the parameter choice rule and render it essentially useless (the simple average of all approximation errors is better than the error err$_{k^*}$ of our algorithm). The left plot, on the other hand, shows that small deviations have hardly any influence. However, in order to illustrate this sensitive dependence on the polynomial exactness of the quadrature rule, we switched from the EGM2008 gravity potential model to potentials with Fourier coefficients that are generated randomly (i.e., in the mean, the Fourier coefficients are equally large at all spherical harmonic degrees). Figure \ref{fig:test4} shows that for a more realistic scenario like the EGM2008 model, the influence of the polynomial exactness of the quadrature rule is significantly smaller. In order to detect a severe failure of our algorithm, we had to decrease the polynomial exactness to degree 80 (opposed to degree 280, which would guarantee the required condition \eqref{cond_3}). This stability of the algorithm is due to the fact the the Fourier coefficients of the EGM2008 gravity potential show a fast decay for growing spherical harmonic degrees. The generally higher optimal errors err$_{\textnormal{opt}}$ in Figures \ref{fig:test3}, \ref{fig:test4} compared to Figures \ref{fig:test1}, \ref{fig:test2} have to be accounted to the influence of the decreased accuracy of the quadrature rule on the approximations $u_k$ via \eqref{eqn:uk} but not to the parameter choice method presented in this paper.

{\begin{figure}
 \scalebox{0.43}{\includegraphics*{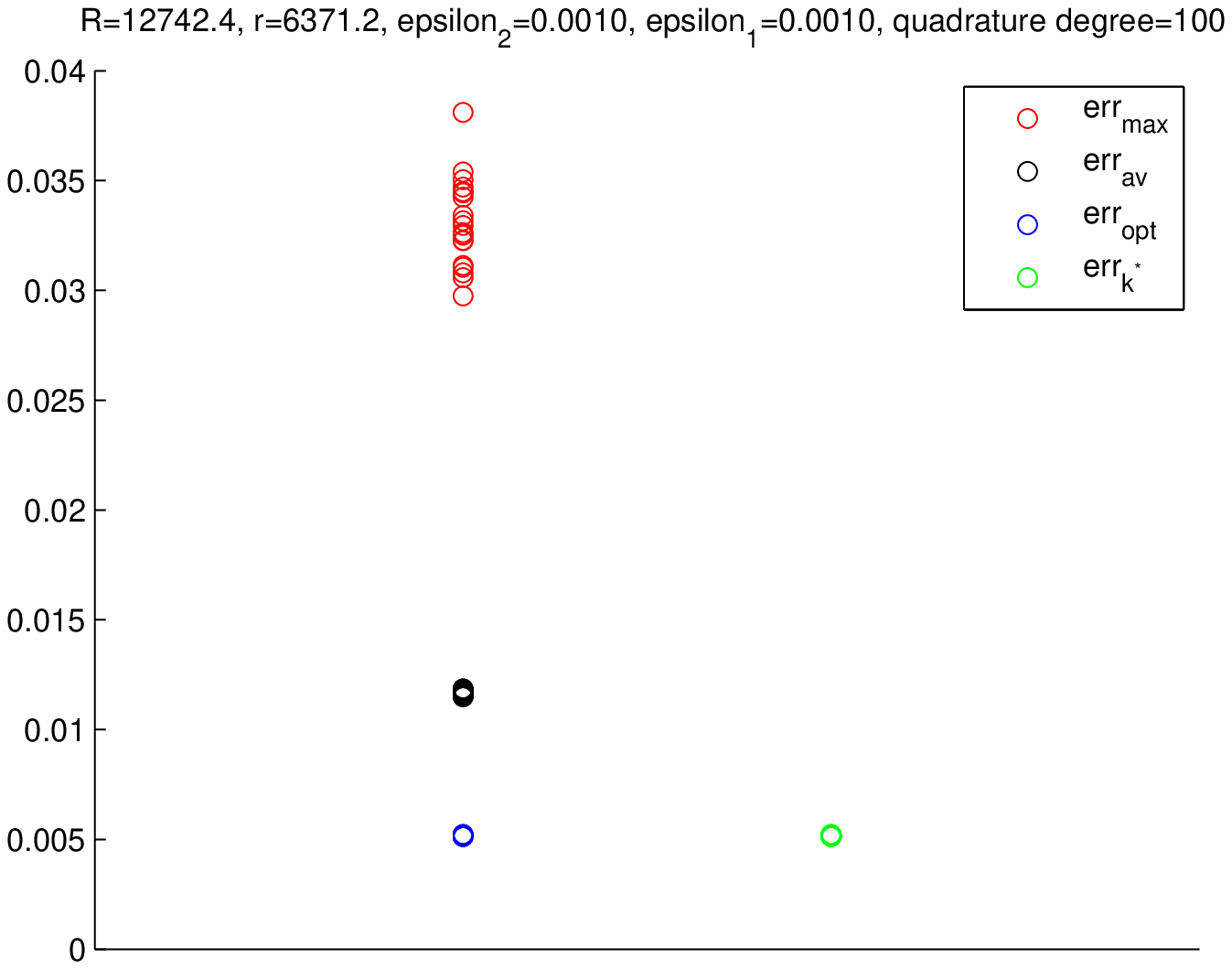}} \scalebox{0.43}{\includegraphics*{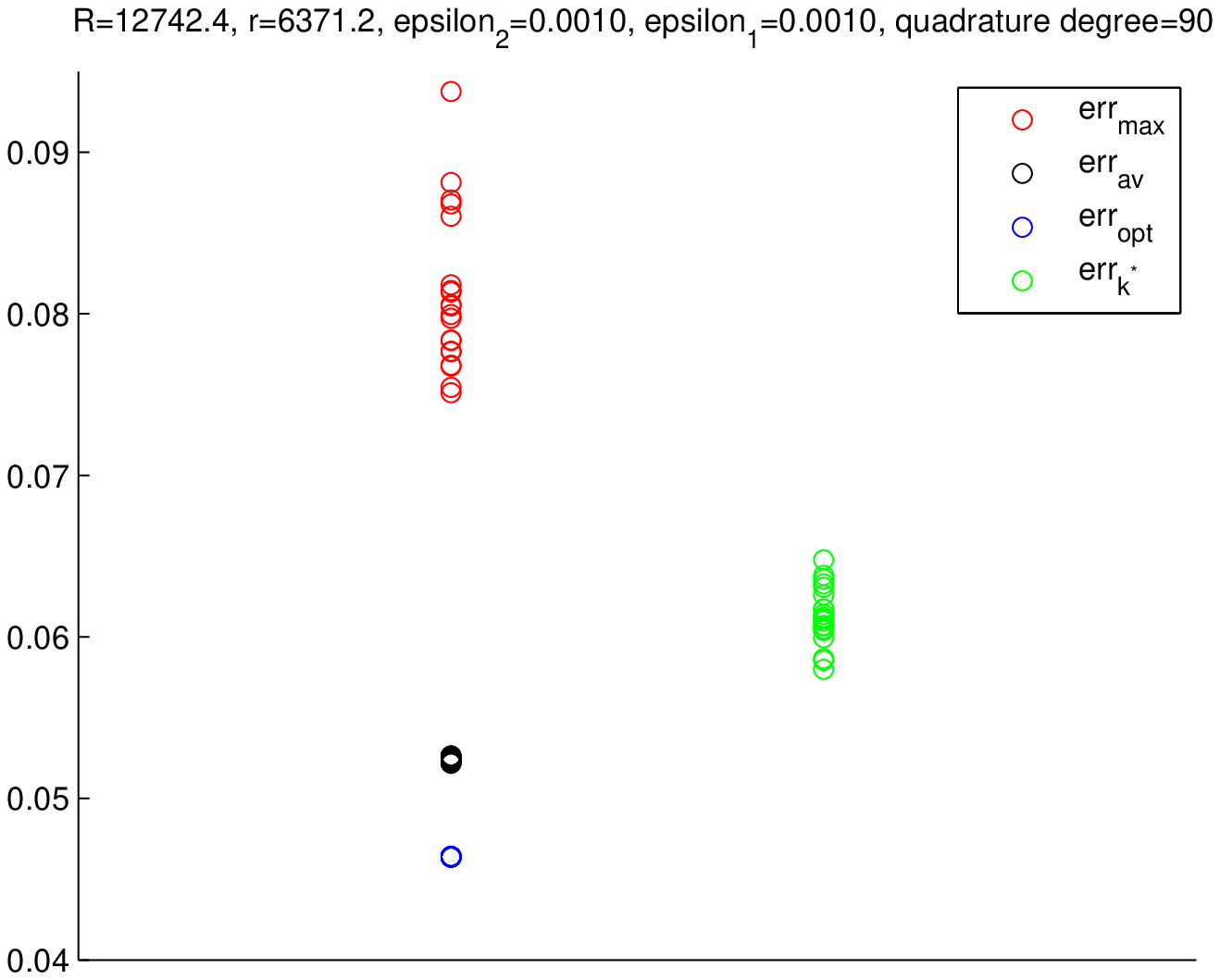}} 
\caption{Relative Errors for the low spherical harmonic degree tests (Situation (a')(1)) for $\varepsilon_1=\varepsilon_2=0.001$ and a quadrature rule with polynomial exactness of degree $100$ (left) and of degree $90$ (right).}\label{fig:test3}
\scalebox{0.43}{\includegraphics*{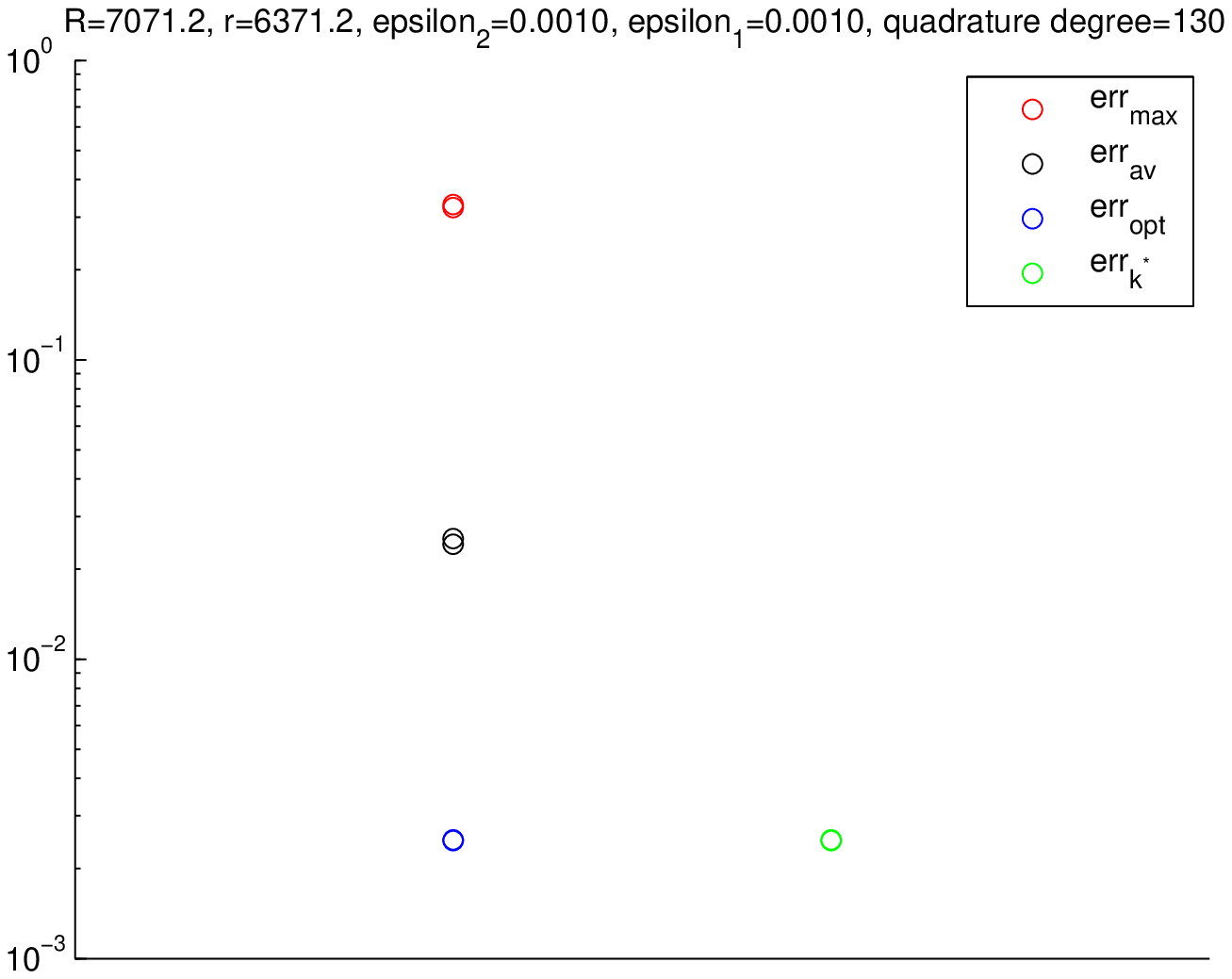}} \scalebox{0.43}{\includegraphics*{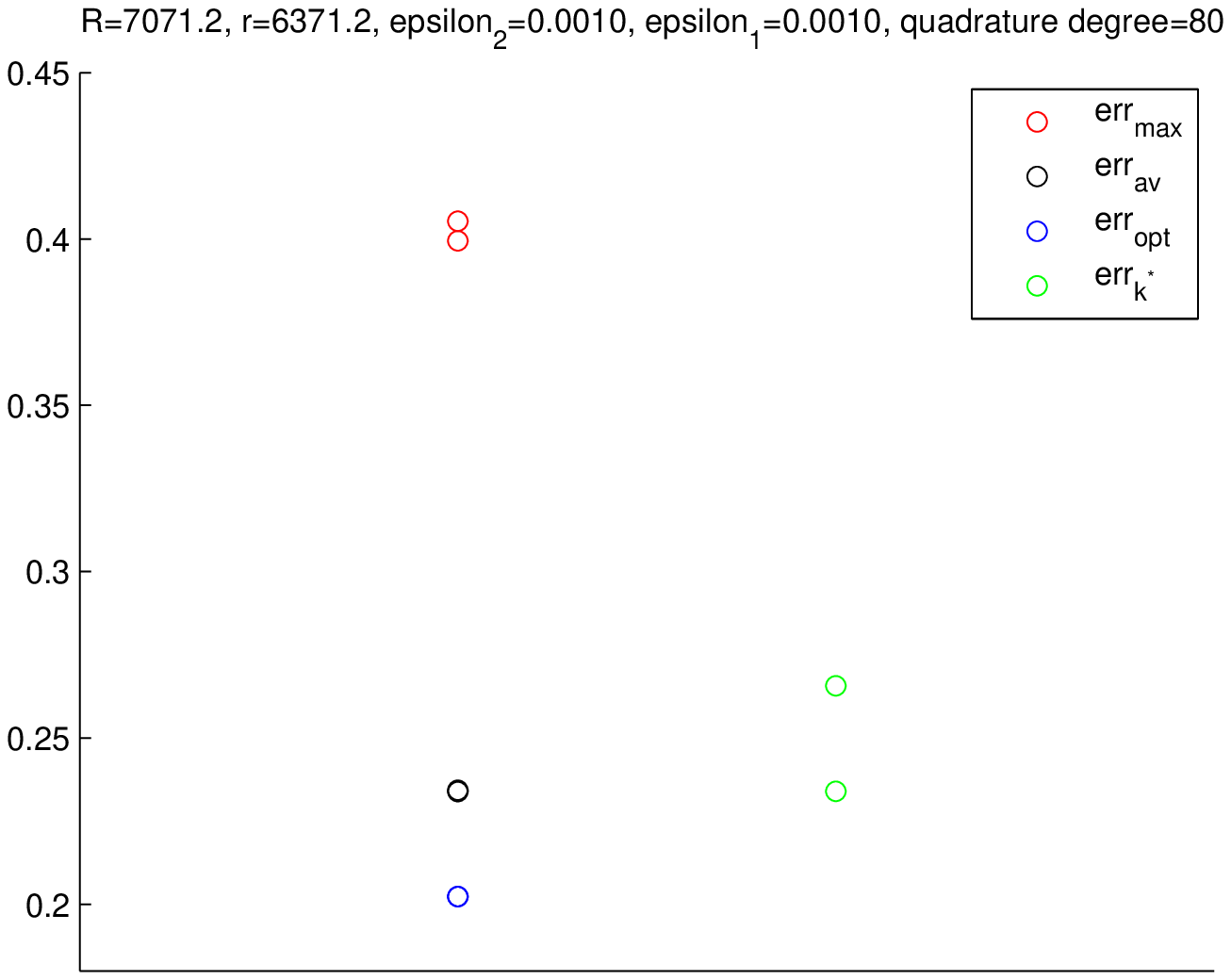}} 
\caption{Relative Errors for the high spherical harmonic degree tests (Situation (a')(2)) for $\varepsilon_1=\varepsilon_2=0.001$ and a quadrature rule with polynomial exactness of degree $130$ (left) and of degree $80$ (right).}\label{fig:test4}
 \end{figure}

\section{Conclusion}
We introduced a simple method to choose a 'good' candidate $u_{k^*}$ among a set of approximations $\{u_k\}_{k=1,\ldots,N}$ of $u^\dag$
and supplied some error estimates for $u_{k^*}$ in relation to $u_{k_{opt}}$.
The numerical illustrations show its good performance and stability for applications, e.g., to the Earth's gravity potential.

\section*{Acknowledgements}

Pavlo Tkachenko gratefully acknowledges the support of the Austrian Science Fund (FWF): project P25424.



\end{document}